\documentclass[reqno]{amsart}
\usepackage{enumitem}

\usepackage[all]{xy}

\usepackage{graphicx}
\usepackage{tikz}
\usetikzlibrary{intersections}
\usetikzlibrary {arrows.meta}
\usetikzlibrary{calc}
\usetikzlibrary{spath3}
\usetikzlibrary{decorations.markings}
\usetikzlibrary{decorations.pathreplacing}
\tikzset{
    position label/.style={
       text height = 1.5ex,
       text depth = 1ex
    },
   brace/.style={
     decoration={brace, mirror},
     decorate
   }
}
\usepackage{mathrsfs}
\usepackage{mathtools}
\usepackage{amsmath}
\usepackage{amssymb}
\usepackage{caption}
\usepackage{subcaption}
\usepackage{hyperref}

\hypersetup{%
  colorlinks=true,
  linkcolor=blue
}

\mathtoolsset{showonlyrefs=true}
\numberwithin{equation}{section}

\newtheorem{thm}{Theorem}[section]
\newtheorem{lem}[thm]{Lemma}
\newtheorem{pro}[thm]{Proposition}

\newtheorem{claim}{Claim}

\theoremstyle{definition}

\newtheorem{rk}[thm]{Remark}

\def\R{{\mathbb R}}
\def\Z{{\mathbb Z}}

\def\bb{\begin}
\def\be{\begin{equation}}
\def\ee{\end{equation}}
\def\bea{\begin{eqnarray}}
\def\eea{\end{eqnarray}}
\def\beaa{\begin{eqnarray*}}
\def\eeaa{\end{eqnarray*}}

\def\bb{\begin}
           \def\ea{\end{array}}
          \def\ec{\end{center}}
     \def\ed{\end{description}}
\def\be{\bb{equation}}        \def\ee{\end{equation}}
\def\bea{\bb{eqnarray}}       \def\eea{\end{eqnarray}}
\def\beaa{\bb{eqnarray*}}     \def\eeaa{\end{eqnarray*}}
 \def\et{\end{thebibliography}}

      \def\e{\varepsilon}

\def\qed {\hfill $\Box$\vskip5pt}

\def\CR{{\rm CR}}

\def\dist{{\rm dist}}

\def\Rep{{\rm Rep}}

\def\Ps{{\rm Ps}}

\newcommand{\capitalize}[1]{\expandafter\capitalizetwo#1\relax}
\def\capitalizetwo#1#2\relax{\MakeUppercase{#1}\MakeLowercase{#2}}
\title{A shadowable chain recurrent set with an attached hyperbolic singularity}

\author[Sogo Murakami]{Sogo Murakami\\
Graduate School of Mathematical Sciences, University of Tokyo, Japan}

\address{Graduate School of Mathematical Sciences, University of Tokyo, Japan}
\email{murakami-sogo880@g.ecc.u-tokyo.ac.jp}



\begin{document}

\begin{abstract}
  We prove that every factor map between topological flows preserves the standard shadowing property if it is injective except for a closed orbit that shrinks to a singularity.
  As an application, we construct a $C^\infty$-flow on a four-dimensional sphere whose nonwandering set contains an attached hyperbolic singularity yet possesses the standard shadowing property. 
  This gives a counterexample to a conjecture given by Arbieto, L\'{o}pez, Rego and S\'{a}nchez (Math. Annalen 390:417-437).
\end{abstract}

\maketitle{}

\section{Introduction}
The shadowing property is a fundamental concept in the theory of dynamical systems, providing a bridge between approximate trajectories and true ones. Informally speaking, for certain types of dynamical systems, any approximate trajectory (commonly referred to as a pseudotrajectory) can be closely followed, or shadowed, by a true orbit.

This property has profound applications across areas: for instance, in the numerical analysis, it is related to sensitivity analysis \cite{WANG2014210}; and in the ergodic theory, it has an application in iterated function systems \cite{Wu:2018wo}. These examples highlight its theoretical and practical versatility in analyzing complex systems.

The shadowing property is closely linked to the structural stability or hyperbolicity. 
Many results concerning the relationship between various shadowing and the stability has been obtained for flows.
Notably, Palmer, Pilyugin and Tikhomirov proved that the Lipschitz shadowing property and the structural stability are equivalent for $C^1$ vector fields on a compact manifold \cite{P.S.T.}. 

On the other hand, beyond the context of uniformly hyperbolic systems, checking the shadowability is an important problem.
For example, the geometric Lorenz attractor does not exhibit the shadowing property as proved by Komuro \cite{KomuroLorenz}. This suggests the limitations of the shadowing in the non-uniformly hyperbolic systems. Indeed, Wen and Wen \cite{WenWenNoShadowing} showed that every singular hyperbolic chain recurrent set with a singularity does not admit the shadowing property. These results showed that the shadowing property is not common in dynamical behavior with singularities beyond the uniform hyperbolicity.

One reason for the lack of the shadowing property is thought to be the presence of singularities.
In this context, Arbieto, L\'{o}pez, Rego and S\'{a}nchez \cite{ArbietoAttachedSing} demonstrated that chain recurrent sets with attached hyperbolic singularities satisfying a certain condition fail to exhibit the shadowing property. Moreover, they conjectured that any chain recurrent set with attached hyperbolic singularities cannot possess this property, suggesting a potential difficulty of possessing the shadowing property in the dynamical systems with singularities.

In this paper, we construct a $C^\infty$ vector field on a four-dimensional sphere whose nonwandering set has an attached singularity, simultaneously satisfies the shadowing property. As a consequence, we disprove the conjecture mentioned above.
The idea of the proof is to construct a ``factor'' of the suspension flow of the Smale horseshoe, exhibiting properties distinct from those observed in systems beyond uniform hyperbolicity so far.

Let $M$ be a $C^\infty$ Riemannian closed manifold with the metric $\dist$ induced by the Riemannian metric.
Let $\phi$ be a $C^\infty$ flow on $M$.

We say that $\xi  : \R \to M$ is a {\it $d$-pseudotrajectory} of $\phi$ if
\[
  \dist \bigl( \xi(t + s), \, \phi(s, \xi(t)) \bigr) < d
\]
for all $t \in \R$ and $s \in [0, 1]$.
Let $\Ps(d)$ be the set of all $d$-pseudotrajectories of $\phi$.
Denote by $\Rep$ the set of all homeomorphisms from $\R$ to $\R$ preserving the orientation.
For $\e > 0$, let
\[
  \Rep(\e) = \left\{ f \in \Rep ; \left\lvert \frac{f(a) - f(b)}{a - b} - 1 \right\rvert < \e, \forall a, b \in \R, a > b \right\}.
\]
We say that a flow $\phi$ has the {\it standard shadowing property} on a $\phi$ invariant set $\Lambda$ if for every $\e > 0$ there exists $d > 0$ such that
if $\xi \in \Ps(d)$ satisfies $\xi(t) \in \Lambda$ for all $t \in \R$ then
\[
  \dist \bigl( \xi(t), \, \phi(h(t), x) \bigr) < \e, \quad t \in \mathbb{R},
\]
for some $x \in \Lambda$ and $h \in \Rep(\e)$.
As a weaker form of the shadowing properties, we say that a flow $\phi$ has the {\it oriented shadowing property} on a $\phi$-invariant set $\Lambda$ if for every $\e > 0$ there exists $d > 0$ such that
if $\xi \in \Ps(d)$ satisfies $\xi(t) \in \Lambda$ for all $t \in \R$ then
\[
  \dist \bigl( \xi(t), \, \phi(h(t), x) \bigr) < \e, \quad t \in \mathbb{R},
\]
for some $x \in \Lambda$ and $h \in \Rep$.

A point $x \in M$ is called a {\it chain recurrent point} of a flow $\phi$ if for any $d, T > 0$ there exists a $d$-pseudotrajectory $\xi$ such that $\xi(0) = \xi(t) = x$ for some $t \geq T$.
Let $\CR(\phi)$ be the set of all chain recurrent points of $\phi$, which is called the {\it chain recurrent set} of $\phi$.
Denote by $\Omega(\phi)$ the nonwandering set of $\phi$. Note that $\Omega(\phi) \subset \CR(\phi)$.
A singularity $p$ (i.e., a fixed point for the flow $\phi$) is {\it attached to} a $\phi$-invariant set $\Lambda$ if $p$ is accumulated by regular points of $\Lambda$, where a regular point is a nonsingular point.
\begin{thm}\label{thm.main}
  There is a $C^\infty$ flow $\phi$ on $S^4$ such that a hyperbolic singularity is attached to $\CR(\phi)$ and $\phi$ has the standard shadowing property on $\CR(\phi)$. Moreover, the chain recurrent set $\CR(\phi)$ coincides with the nonwandering set $\Omega(\phi)$ of $\phi$.
\end{thm}
\begin{rk}
  Theorem \ref{thm.main} gives a counterexample to Conjecture $2$ by Arbieto, L\'{o}pez, Rego and S\'{a}nchez \cite{ArbietoAttachedSing}.
  They conjectured that a chain recurrent set with an attached hyperbolic singularity does not have the oriented shadowing property.
  Since the standard shadowing property implies the oriented shadowing property, our result is stronger than the conjecture.
  In addition, the answer to Question $1$ of \cite{ArbietoAttachedSing} can be also constructed as a topological flow on a closed disk. See remark in \cite[Section 1]{Mura}.
\end{rk}
The chain recurrent set $\CR(\phi)$ of Theorem \ref{thm.main} is a ``factor'' of the chain recurrent set of a modified suspension flow of Smale horseshoe.
Theorem \ref{thm.main} is proved via the following theorem:
\begin{thm}\label{thm.factor}
  Consider topological flows $(\Lambda_1, \phi_1)$ and $(\Lambda_2, \phi_2)$ on compact metric spaces.
  Suppose that there is a continuous surjection $P : \Lambda_1 \to \Lambda_2$ by which the following diagram commutes:
  \[
    \xymatrix{
    \Lambda_1 \ar[r]^{\phi_1(t,\cdot)} \ar[d]_{P} & \Lambda_1 \ar[d]^{P} \\
    \Lambda_2 \ar[r]_{\phi_2(t,\cdot)} & \Lambda_2 \ar@{}[lu]|{\circlearrowright}
    }
  \]
  for all $t \in \R$.
  If there exists a closed orbit $\gamma \subset \Lambda_1$ of $\phi_1$ such that $P(\gamma)$ is a singularity of $\phi_2$, and the restriction of $P$ to $\Lambda_1 \setminus \gamma$ is a bijection onto $\Lambda_2 \setminus P(\gamma)$, then the standard shadowing property of $\phi_1$ implies that of $\phi_2$.
\end{thm}

\section{Preliminaries}\label{sec.preliminaries}
\begin{figure}[h]
  \includegraphics[width=12cm]{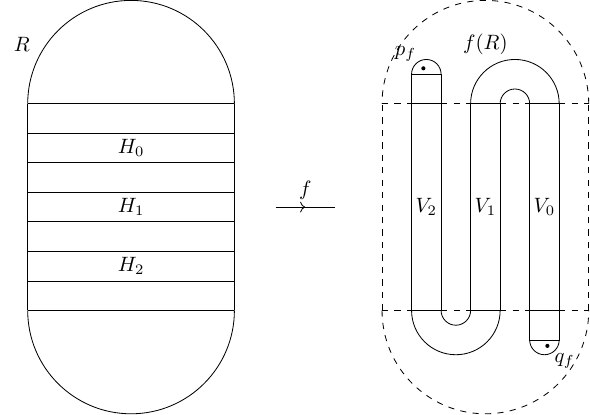}
  \caption{$3$-fold Smale horseshoe with $f(H_i) = V_i$ for $i = 0,1,2$.}\label{figure.Horseshoe}
\end{figure}

In this section, we describe the structure of the desired flow on $S^4$ and reduce the proof of Theorem \ref{thm.main} to that of Theorem \ref{thm.factor}.
Let us denote by $\| \cdot \|$ the Euclidean norm of the Euclidean space.
Now, we consider a $3$-fold Smale horseshoe diffeomorphism $f : \R^2 \to \R^2$, which will play an important role in the construction of our flow.
Take a shape $R$ that consists of a unit square centered at the origin of $\R^2$ with semicircles attached to its top and bottom sides, then we may assume $f(R) \subset R$ (see Figure \ref{figure.Horseshoe}).
The horizontal strips $H_0$, $H_1$, $H_2$ are mapped linearly onto vertical strips $V_0$, $V_1$, $V_2$, respectively.
We may also assume that:
\begin{itemize}
  \item The chain recurrent set $\CR(f)$ of $f$ consists of an attracting periodic orbit $\{ p_f, q_f \}$ of period $2$ (see Figure \ref{figure.Horseshoe}) and a transitive non-trivial hyperbolic set $\Lambda_f$.
  \item $f(x,y) = (-x/2,-y/2)$ for all $(x,y) \in \R^2$ with sufficiently large $\| (x,y) \|$.
  \item For all $(x,y) \in \R^2$, we have $f(-x,-y) = -f(x,y)$.
  \item For all $(x,y) \in \R^2$, there exists $n \geq 1$ such that $f^n(x,y) \in R$.
\end{itemize}

Let $\bar{f} : \R^2 \to \R^2$ be
\begin{equation}
  \bar{f}(x,y) = -f(x,y)\label{eq.def_of_ol_f}
\end{equation}
and let $\tilde{\phi} : \R \times \widetilde{M} \to \widetilde{M}$ be a modified suspension flow of $\bar{f}$ (see Section \ref{sec.CRset_of_phi} for the precise definitions of $\tilde{\phi}$ and $\widetilde{M}$).
In Section \ref{sec.constructionofX}, we will construct a $C^\infty$ vector field $X$ on $\R^4$ satisfying the following properties:
\begin{enumerate}
  \item[\bf{(P1)}] $(0,0,0,0) \in \R^4$ is the only singularity of $X$, which is a hyperbolic singularity of index $2$.
  \item[\bf{(P2)}] Let $\phi : \R \times \R^4 \to \R^4$ be the flow generated by $X$.
        There exists a smooth map $\widetilde{P} : \widetilde{M} \to \R^4$ such that the following diagram commutes:
        \begin{equation}
          \xymatrix{
          \widetilde{M} \ar[r]^{\tilde{\phi}(t,\cdot)} \ar[d]_{\widetilde{P}} & \widetilde{M} \ar[d]^{\widetilde{P}} \\
          \R^4 \ar[r]_{\phi(t,\cdot)} & \R^4 \ar@{}[lu]|{\circlearrowright}
          }\label{eq.cumm_diag_phi}
        \end{equation}
        for all $t \in \R$.
  \item[\bf{(P3)}] $\widetilde{P}(\CR(\tilde{\phi})) = \CR(\phi)$.
  \item[\bf{(P4)}] $\widetilde{P}^{-1}(0,0,0,0)$ is a hyperbolic closed orbit of $\tilde{\phi}$ and is accumulated by regular points of $\CR(\tilde{\phi})$.
        On the other hand, the restriction of $\widetilde{P}$ to the domain $\widetilde{M} \setminus \widetilde{P}^{-1}(0,0,0,0)$ is injective.
  \item[\bf{(P5)}] There is $N_0 > 0$ such that
        \[
          X(x,y,z,w) = (-x-\pi y, -y+\pi x, -z-\pi w, -w+\pi z)
        \]
        for all $(x,y,z,w) \in \R^4$ with $\| (x,y,z,w) \| \geq N_0$.
\end{enumerate}

We construct a vector field $X_{S^4}$ on $S^4$ from $X$.
Let $F : \R^4 \to S^4$ be
\[
  F(p) =
  \left(
  \frac{4x}{\| p \|^2+4},
  \frac{4y}{\| p \|^2+4},
  \frac{4z}{\| p \|^2+4},
  \frac{4w}{\| p \|^2+4},
  \frac{\| p \|^2-4}{\| p \|^2+4}
  \right),
\]
where $p = (x,y,z,w) \in \R^4$.
Then $F$ is a $C^\infty$ diffeomorphism from $\R^4$ to $S^4 \setminus \{(0,0,0,0,1)\}$.
Define a vector field $X_{S^4}$ on $S^4$ by 
\[
  X_{S^4}(p) = 
  \begin{cases}
    F_*(X)(p), & p \in F(\R^4),    \\
    0,         & p = (0,0,0,0,1).
  \end{cases}
\]
Note that, if $X_{S^4}$ has the shadowing property on its chain recurrent set $\CR(X_{S^4})$, then $\CR(X_{S^4}) = \Omega(X_{S^4})$. 
Therefore, it suffices to show that $\CR(X_{S^4})$ has an attached hyperbolic singularity and $X_{S^4}$ has the shadowing property on $\CR(X_{S^4})$.
It is easy to see that $(0,0,0,0,1)$ is a hyperbolic repeller of $X_{S^4}$. In fact, the real parts of all eigenvalues of $DX_{S^4}$ at $(0,0,0,0,1)$ are positive from property {\bf (P5)}.
Thus, $\CR(X_{S^4}) = F(\CR(X)) \cup \{(0,0,0,0,1)\}$.
Since it is obvious that $X_{S^4}$ has the standard shadowing property on a hyperbolic repeller,
the proof of Theorem \ref{thm.main} is reduced to showing the following proposition:
\begin{pro}
  $F(\CR(X))$ has an attached hyperbolic singularity and $X_{S^4}$ has the standard shadowing property on $F(\CR(X))$.
\end{pro}
By the definition of $X_{S^4}$, this proposition follows from the following proposition:
\begin{pro}\label{pro.reduction_to_R4}
  $\CR(X)$ has an attached hyperbolic singularity and $\phi$ has the standard shadowing property on $\CR(X)$.
\end{pro}
By properties {\bf(P2)}, {\bf(P3)} and {\bf(P4)}, we can apply Theorem \ref{thm.factor} to $\widetilde{P} \vert_{\CR(\tilde{\phi})} : \CR(\tilde{\phi}) \to \CR(\phi)$ as $\Lambda_1 = \CR(\tilde{\phi})$ and $\Lambda_2 = \CR(\phi)$ in order to prove that $\phi$ has the standard shadowing property on $\CR(X)$ (note that the suspension flow of a diffeomorphism with the standard shadowing property also has the standard shadowing property by \cite[Theorem 2]{TomasSuspensionShadowing}).
It follows from property {\bf(P4)} that $(0,0,0,0)$ is an attached hyperbolic singularity, and thus we prove Proposition \ref{pro.reduction_to_R4}, which finishes the proof of Theorem \ref{thm.main}.

The structure of the paper is as follows.
In Section \ref{sec.constructionofX}, we construct the vector field $X$ and check that properties {\bf(P1)} and {\bf(P5)} hold,
and in Section \ref{sec.CRset_of_phi}, we construct $\widetilde{M}$ and $\tilde{\phi}$ and prove that $\CR(X)$ is a compact invariant set satisfying properties {\bf(P2)}, {\bf(P3)} and {\bf(P4)} above.
In Section \ref{sec.prf_of_thm_factor}, we prove Theorem \ref{thm.factor}.

\section{Construction of $X$}\label{sec.constructionofX}
\begin{figure}[h]
 \includegraphics[width=12cm]{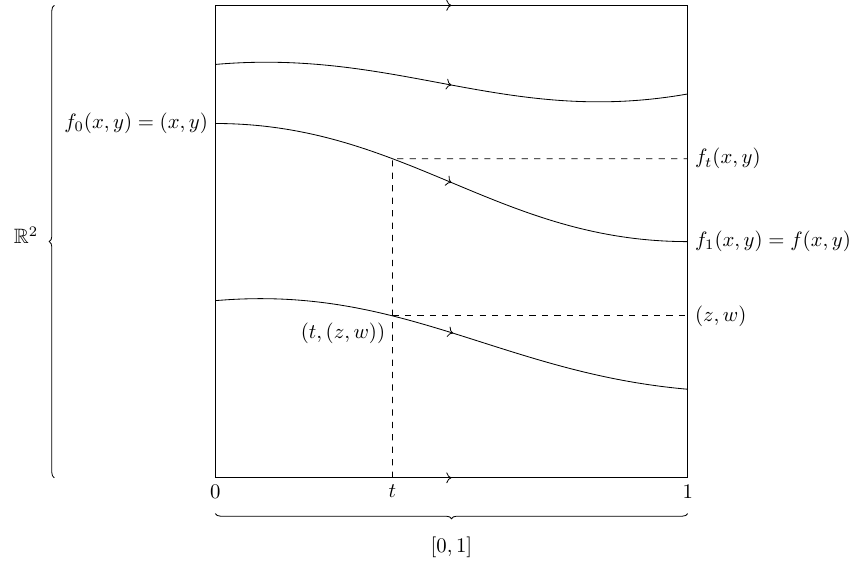}
  \caption{A family of diffeomorphisms $\{f_t\}_{0 \leq t \leq 1}$. The slope of the curve at $(t, (z, w))$ is $V(z, w, t)$.}\label{figure.suspensionoff}
\end{figure}

In this section, we construct a vector field $X$ satisfying properties {\bf(P1)} and {\bf(P5)}.
Let us consider a family of diffeomorphisms $\{f_t : \R^2 \to \R^2\}_{0 \leq t \leq 1}$ satisfying the following properties (see Figure \ref{figure.suspensionoff}):
\begin{enumerate}
  \item $\{f_t : \R^2 \to \R^2\}_{0 \leq t \leq 1}$ is a smooth transformation from ${\rm id}_{\R^2}$ to $f$ given at the beginning of Section \ref{sec.preliminaries}; i.e., $f_0 \equiv {\rm id}_{\R^2}$, $f_1 \equiv f$ and the map
        \[
          \R^2 \times \R/\Z \ni (p,t) \mapsto
          \frac{\partial f_t}{\partial t}(p) \in T_{f_t(p)}\R^2
        \]
        is smooth.
  \item $f_t(x,y) = -f_t(-x,-y)$.
  \item There is $\nu_0 > 0$ such that
        \[
          \frac{\partial f_t}{\partial t}(x,y) = (-(\log 7) x,(\log 7) y)
        \]
        for all $(x,y) \in \R^2$ and $t \in [0,1]$ with $\| f_t(x,y) \| \leq \nu_0$.
  \item There is $N_0 > 0$ such that
        \[
          \frac{\partial f_t}{\partial t}(x,y) = (-x,-y)
        \]
        for all $(x,y) \in \R^2$ and $t \in [0,1]$ with $\| f_t(x,y) \| \geq N_0/2$.
\end{enumerate}
Let $\bar{f}$ be the map defined in \eqref{eq.def_of_ol_f}.
Using $\{f_t\}$, we may construct a modified suspension flow of $\bar{f}$.
For $(x,y) \in \R^2$ and $t \in [0,1]$, define
\[
  V(x,y,t) = \frac{\partial}{\partial s}\bigg|_{s=0} f_{t+s} \circ f_t^{-1}(x,y) \in \R^2.
\]
This corresponds to the slope of the flow at the point $(t,(x,y))$ in Figure \ref{figure.suspensionoff}.

Let us construct a vector field $X$ on $\R^4$ satisfying the conditions given at the end of Section \ref{sec.preliminaries}.
For $\theta \in [0,2]$, define $e(\theta) = (\cos \pi\theta,\sin \pi\theta) \in \R^2$.
In this paper, we regard a point in $\R^4$ as a point in $(\R^2)^2$ to express that point by a pair of polar coordinates.
For $p,q \in \R$ and $\theta,\varphi \in [0,2]$, let
\[
  R_{\theta, \varphi}(p,q) = (p \cdot e(\theta),\, q \cdot e(\varphi)) \in \R^4.
\]
The vector field $X$ on $\R^4$ is defined by the linear sum of three vector fields $Y$, $Z$, $W$.
We define these vector fields separately as follows.

Let $Y(x,y,z,w) = (-\pi y,\, \pi x,\, -\pi w,\, \pi z)$.
Notice that $Y(a \cdot e(\theta),\, b \cdot e(\varphi)) = (\pi a \cdot e(\theta+1/2),\, \pi b \cdot e(\varphi+1/2))$.

Letting
\begin{align}
  D 
   & = \{ R_{\theta,\theta}(a,b) ; a,b \in \R, \,\theta \in [0,2] \}                                                  \\
   & \cup \{ R_{\theta,\varphi}(a,b) ; a,b \in \R, \, \theta \in [0,2], \| (a,b) \| \leq \nu_0/2 \}                   \\
   & \cup \{ R_{\theta,\varphi}(a,b) ; a,b \in \R, \, \theta \in [0,2], \| (a,b) \| \geq N_0 \},\label{eq.def_of_DZ0}
\end{align}
we define a smooth vector field $Z_0$ on $D$ by
\[
  Z_0(R_{\theta,\varphi}(a,b)) =
  \begin{cases}
    R_{\theta,\varphi}(V_1(a,b,\theta),V_2(a,b,\theta)), & \theta = \varphi,         \\
    R_{\theta,\varphi}(-(\log 7)a,(\log 7)b),            & \| (a,b) \| \leq \nu_0/2, \\
    R_{\theta,\varphi}(-a,-b),                           & \| (a,b) \| \geq N_0
  \end{cases}
\]
where $V_1, V_2 : \R^2 \times [0,2] \to \R$ are functions such that $V(a,b,\theta) = (V_1(a,b,\theta),V_2(a,b,\theta))$.
Then by applying \cite[Lemma 2.26]{lee2012smooth} to $Z_0$, we obtain a smooth vector field $Z$ on $\R^4$ satisfying $Z(p) = Z_0(p)$ for all $p \in D$.


Let us define a vector field $W_0$ on $\R^4$ by
\begin{align}
  W_0(x,y,z,w)
   & = -(xw-yz) \, {\rm grad}(xw-yz)                          \\
   & = \left( -(xw-yz)w,(xw-yz)z,(xw-yz)y,-(xw-yz)x \right).
\end{align}
To localize this vector field, we introduce a smooth bump function $\tau_W : \R_{\geq 0} \to [0,1]$ such that
\[
  \tau_W(r) =
  \begin{cases}
    1, & r \leq N_0/2, \\
    0, & r \geq N_0.
  \end{cases}
\]
Then define a vector field $W$ by
\[
  W = \tau_W W_0.
\]
Next, define $h : \R^4 \to \R$ by
\[
  h(x,y,z,w) = xw-yz.
\]
The following lemma asserts that every integral curve of $Y$ is contained in a level set of $h$ (for the definition of a level set, see \cite[Example 1.32]{lee2012smooth}).
\begin{lem}\label{lem.boundforY}
  For all $p = (x,y,z,w) \in \R^4$,
  \[
    Dh(p) (Y(p)) = 0.
  \]
\end{lem}
\begin{prf}
  Denote by $\langle \cdot, \cdot \rangle_{\R^4}$ the inner product of elements of $\R^4$.
  Then
  \begin{align}
    Dh(p) (Y(p))
     & = \left\langle 
    \left(\frac{\partial h}{\partial x}(p),
    \frac{\partial h}{\partial y}(p),
    \frac{\partial h}{\partial z}(p),
    \frac{\partial h}{\partial w}(p) \right),
    Y(p) \right\rangle_{\R^4}                           \\
     & = \left\langle 
    \left(w,
    -z,
    -y,
    x \right),
    Y(p) \right\rangle_{\R^4}                           \\
     & = w(-\pi y) + (-z)\pi x + (-y) (-\pi w) + x\pi z \\
     & = 0.
  \end{align}
  \qed
\end{prf}
The following two lemmas gives an evaluation of the change of $h$ along the integral curves of $Z$ and $W$, respectively.
\begin{lem}\label{lem.boundforZ}
  There exists $C > 0$ such that
  \[
    \lvert Dh(p) (Z(p)) \rvert \leq C \lvert h(p) \rvert
  \]
  for all $p \in \R^4$.
  Moreover, 
  \[
    Dh(p) (Z(p)) = 0
  \]
  for all $p \in \R^4$ with $\| p \| \leq \nu_0/2$.
\end{lem}
\begin{prf}
  Notice that for all $p = (x,y,z,w) = (a \cdot e(\theta), b \cdot e(\varphi)) \in \R^4$, it follows that $Dh(p) = (w,-z,-y,x) = (b \cdot e(\varphi-1/2), a \cdot e(\theta+1/2))$.
  We have
  \begin{align}
    Dh(p) (Z(p)) 
     & = Dh(p) (Z_0(p))                                                              \\
     & = \langle (b \cdot e(\varphi-1/2),\, a \cdot e(\theta+1/2)),\, Z_0(p) \rangle \\
     & =
    \begin{cases}
      0,      & \theta = \varphi,         \\
      0,      & \| (a,b) \| \leq \nu_0/2, \\
      -2h(p), & \| (a,b) \| \geq N_0      \\
    \end{cases}\label{eq.DhZ0}
  \end{align}
  for all $p \in D$, where $D$ is the set given in \eqref{eq.def_of_DZ0}.
  Thus, the proof of the lemma is reduced to proving the following claim.
  \begin{claim}
    There is $C \geq 2$ such that
    \[
      \lvert Dh(p)(Z(p)) \rvert \leq C \lvert h(p) \rvert
    \]
    for all $p \in \overline{\R^4 \setminus D}$.
  \end{claim}
  Assume to the contrary that, there are $2 \leq C_1 < C_2 < \cdots \to \infty$ and $p_n \in \overline{\R^4 \setminus D}$ such that $\lvert Dh(p_n)(Z(p_n)) \rvert > C_n \lvert h(p_n) \rvert$ for all $n \geq 1$.
  By the compactness, there is $p_0 \in \overline{\R^4 \setminus D}$ satisfying $p_{n_k} \to p_0$ as $k \to \infty$.
  For all $k$ and $l$ with $k > l$, we have
  \[
    \lvert Dh(p_{n_k})(Z(p_{n_k})) \rvert
    > C_{n_k} \lvert h(p_{n_k}) \rvert
    \geq C_{n_{l}} \lvert h(p_{n_k}) \rvert.
  \]
  Taking $k \to \infty$, we obtain
  \[
    \lvert Dh(p_0)(Z(p_0)) \rvert \geq C_{n_{l}} \lvert h(p_0) \rvert.
  \]
  Consequently, we have $\lvert h(p_0) \rvert = 0$.
  If $p = (a \cdot e(\theta), b \cdot e(\varphi)) \in \R^4$ satisfies $h(p) = 0$, then $\theta = \varphi$.
  From \eqref{eq.DhZ0}, we have $Dh(p)(Z(p)) = 0$ for all $p \in \R^4$ with $\theta = \varphi$. 
  Thus, if $h(p) = 0$ then $Dh(p)(Z(p)) = 0$.
  Since $p_0 \neq (0,0,0,0)$, the point $p_0$ is a regular point of $h$. 
  Therefore, there exist a neighborhood $U$ of $p$ and a constant $C > 0$ such that $\lvert Dh(p)(Z(p)) \rvert \leq C \lvert h(p) \rvert$ for all $p \in U$.
  This contradicts with the fact that $C_n > C$ for sufficiently large $n$.
  \qed
\end{prf}
\begin{lem}\label{lem.boundforW}
  For all $p \in \R^4$ with $\| p \| \leq N_0/2$,
  \[
    Dh(p) (W(p)) = - h(p) \cdot \| p \|^2.
  \]
\end{lem}
\begin{prf}
  Note that $W(p) = -h(p) \, {\rm grad} \, h(p)$ for all $p \in \R^4$ with $\| p \| \leq N_0/2$.
  Then
  \begin{align}
    Dh(p) (W(p))
     & = \langle (w,-z,-y,x), -h(p) {\rm grad} \, h(p) \rangle_{\R^4} \\
     & = -h(p)\langle (w,-z,-y,x), {\rm grad} \, h(p) \rangle_{\R^4}  \\
     & = -h(p)\langle (w,-z,-y,x), (w,-z,-y,x) \rangle_{\R^4}         \\
     & = -h(p)(x^2+y^2+z^2+w^2)                                       \\
     & = -h(p)\| p \|^2.
  \end{align}
  \qed
\end{prf}
Let $C$ be the constant given in Lemma \ref{lem.boundforZ}.
Define
\[
  X = Y + Z + 2C(\nu_0/2)^{-2}W.
\]
We now explain the roles of the vector fields $Y$, $Z$, and $W$. 
Let $K \subset \R^4$ be the closed set defined at the beginning of Section \ref{sec.CRset_of_phi}.
The vector field $W$ has been introduced to ensure that points outside $K$ are attracted to $K$ (see Proposition \ref{pro.h_is_hightfunction}).
This fact is essential in the proof of property {\bf (P3)}.
Since $W=0$ on $K$, it follows that $X = Y + Z$ on $K$. 
We define $Y$ and $Z$ to satisfy property {\bf (P2)}.

First, we prove that $X$ satisfies properties {\bf (P1)} and {\bf (P5)} in this section.
For every $p = (a \cdot e(\theta), b \cdot e(\varphi)) \in \R^4$ with $\| p \| \leq \nu_0/2$, 
\begin{align}
  X(p)
   & = Y(p) + Z(p) + 2C(\nu_0/2)^{-2}W(p)                                     \\
   & = (-\pi y, \pi x, -\pi w, \pi z)                                         \\
   & + (-\log 7 x, -\log 7 y, \log 7 z, \log 7 w)                             \\
   & + 2C(\nu_0/2)^{-2} \left(-(xw-yz)w,(xw-yz)z,(xw-yz)y,-(xw-yz)x \right).
\end{align}
Thus, we have
\begin{align}
  DX(0,0,0,0)
   & = DY(0,0,0,0) + DZ(0,0,0,0) + 0 \\
   & =
  \begin{pmatrix}
    -\log 7 & -\pi    & 0      & 0      \\
    \pi     & -\log 7 & 0      & 0      \\
    0       & 0       & \log 7 & -\pi   \\
    0       & 0       & \pi    & \log 7 \\
  \end{pmatrix}
\end{align}
and that the eigenvalues of $DX(0,0,0,0)$ is $\pm \log 7 \pm i\pi$, which implies property {\bf (P1)}.

As for property {\bf (P5)},
\begin{align}
  X(p)
   & = Y(p) + Z(p) + 2C(\nu_0/2)^{-2}W(p)                   \\
   & = (-\pi y,\, \pi x,\, -\pi w,\, \pi z)                 \\
   & + (-x, -y, -z, -w)                                     \\
   & + (0,0,0,0)                                            \\
   & = (-x-\pi y,\, -y + \pi x,\, -z -\pi w,\, -w + \pi z)
\end{align}
for all $p = (x, y, z, w) \in \R^4$ with $\| p \| \geq N_0$, proving property {\bf (P5)}.
\section{A chain recurrent set with an attached singularity}\label{sec.CRset_of_phi}
Let $\phi$ be the flow on $\R^4$ generated by $X$.
In this section, we define the modified suspension flow $(\widetilde{M}, \tilde{\phi})$ of $\bar{f}$ (see \eqref{eq.def_of_ol_f} for the definition of $\bar{f}$) and prove properties {\bf (P2)}, {\bf (P3)} and {\bf (P4)}.

Let $K = \{ p \in \R^4 ; h(p) = 0 \}$.
The following proposition implies that every orbit of $\phi$ converges ``monotonically'' to $K$.

\begin{pro}\label{pro.h_is_hightfunction}
  Let ${\rm sgn} : \R \to \{-1, 0, 1\}$ be
  \[
    {\rm sgn}(x) =
    \begin{cases}
      -1, & x < 0,  \\
      0,  & x = 0,  \\
      1,  & x > 0.
    \end{cases}
  \]
  Then
  \[
    {\rm sgn}\left( \frac{\partial}{\partial t} \bigg|_{t=0} h(\phi_t(p)) \right)
    = -{\rm sgn}(h(p))
  \]
  for all $p \in \R^4$ and
  the value $\lvert h(\phi_t(p)) \rvert$ converges monotonically to $0$ as $t \to \infty$.
\end{pro}
\begin{prf}
  Let $C > 0$ be the constant given in Lemma \ref{lem.boundforZ}.
  By Lemmas \ref{lem.boundforY} and \ref{lem.boundforW}, we have
  \begin{align}
    \frac{\partial}{\partial t} \bigg|_{t=0} h(\phi_t(p))
     & = Dh(p) (X(p))                                             \\
     & = Dh(p) (Y+Z+2C(\nu_0/2)^{-2}W)(p)                         \\
     & = Dh(p)Y(p) + Dh(p)Z(p) + 2C(\nu_0/2)^{-2} Dh(p)W(p)       \\
     & = 0 +  Dh(p)Z(p) - 2C(\nu_0/2)^{-2} h(p) \cdot \| p \|^2.
  \end{align}
  On the other hand, by Lemma \ref{lem.boundforZ},
  \[
    \lvert Dh(p)Z(p) \rvert \leq C(\nu_0/2)^{-2} \lvert h(p) \rvert \cdot \| p \|^2
  \]
  regardless of whether $\| p \| \leq \nu_0/2$ holds or not (that is, if $\| p \| \leq \nu_0/2$ then $Dh(p)Z(p) = 0$ and if $\| p \| \geq \nu_0/2$ then $\lvert Dh(p)Z(p) \rvert \leq C\lvert h(p) \rvert \leq C(\nu_0/2)^{-2} \lvert h(p) \rvert \cdot \| p \|^2$).
  Thus,
  \begin{align}
    {\rm sgn}\left( \frac{\partial}{\partial t} \bigg|_{t=0} h(\phi_t(p)) \right)
     & = {\rm sgn}(-C(\nu_0/2)^{-2} h(p) \cdot \| p \|^2) \\
     & = -{\rm sgn}(h(p)).
  \end{align}
  This completes the proof of the proposition.
  \qed
\end{prf}
The following proposition follows from property {\bf (P5)}, whose proof is given at the end of Section \ref{sec.constructionofX}.
\begin{pro}\label{pro.X_attracting_nbd}
  If $p \in \R^4$ satisfies $\| p \| \geq N_0$, then
  \[
    \left\langle X(p), \frac{p}{\| p \|} \right\rangle_{\R^4} = -\| p \|.
  \]
\end{pro}
\begin{prf}
  By property {\bf (P5)}, we have
  \[
    X(x,y,z,w) = (-x-\pi y,-y+\pi x,-z-\pi w,-w+\pi z)
  \]
  for all $p = (x,y,z,w) \in \R^4$ with $\| p \| \geq N_0$.
  Thus, 
  \begin{align}
    \langle X(p), p \rangle_{\R^4}
     & = \langle (-x-\pi y,-y+\pi x,-z-\pi w,-w+\pi z), (x,y,z,w) \rangle_{\R^4} \\
     & = (-x-\pi y)x + (-y+\pi x)y + (-z-\pi w)z + (-w+\pi z)w                     \\
     & = -(x^2+y^2+z^2+w^2),
  \end{align}
  finishing the proof.
  \qed
\end{prf}
\begin{pro}\label{pro.CRphi_bdd_inK}
  The chain recurrent set $\CR(\phi)$ is bounded and $\CR(\phi) \subset K$.
\end{pro}
\begin{prf}
  Let $x_0 \in \CR(\phi)$ and let $\{g_n\}_n$ be $(1/n)$-cycle from $x_0$ (i.e., $g_n(0) = x_0$ and $g_n(t_n) = x_0$ for some $t_n$ with $t_n \to \infty$ as $n \to \infty$).
  It follows from Proposition \ref{pro.X_attracting_nbd} that
  \[
    \phi(1, \{ p \in \R^4 ; \| p \| \leq N_0+1 \})
    \subset \{ p \in \R^4 ; \| p \| \leq N_0 \},
  \]
  that is, $\{ p \in \R^4 ; \| p \| \leq N_0+1 \}$ is an attracting neighborhood of $(0,0,0,0)$.
  Consequently, we can assume that $g_n(\R) \subset \{ p \in \R^4 ; \| p \| \leq N_0 + 1 \}$ for sufficiently large $n$.
  This implies that $x_0 = g_n(0)$ satisfies $\| x_0 \| \leq N_0+1$ and thus $\CR(\phi)$ is bounded.
  
  Assume to the contrary that there is $x_0 \in \CR(\phi) \setminus K$.
  Without loss of generality, we may assume that $h(x_0) > 0$, because the other case where $h(x_0) < 0$ can be proven analogously.
  For $t \geq 0$, let
  \[
    K(t) = \{ p \in \R^4 ; \lvert h(p) \rvert \leq t, \| p \| \leq N_0+1 \}.
  \]
  From Proposition \ref{pro.h_is_hightfunction}, we may choose $\lambda \in (h(\phi(1,x_0)), h(x_0))$.
  Then there exists $N_1$ such that $g_n(1) \in K(\lambda)$ for all $n \geq N_1$.
  Again, using Proposition \ref{pro.h_is_hightfunction},
  we see that $\max \{h(x) ; x \in \phi(1,K(\lambda)) \} < \lambda$. Thus, there exists $d > 0$ such that if $p \in K(\lambda)$ and $q \in \{ p \in \R^4 ; \| p \| \leq N_0+1 \}$ satisfies $\dist(\phi(1,p),q) < d$, then $q \in K(\lambda)$.
  Thus, for $n \geq N_1$ with $1/n < d$ and $k \geq 1$, we have $g_n(k) \in K(\lambda)$.
  Taking $d$ smaller if necessary, this contradicts our choice of $\lambda$ with $g_n(t_n) = x_0 \notin K(\lambda)$.
  \qed
\end{prf}

Now, let us construct the modified suspension flow of Smale horseshoe $(\widetilde{M}, \tilde{f})$, which has been briefly introduced in Section \ref{sec.preliminaries}.
Let us define an equivalence relation
$ (x_1, y_1, t_1) \sim (x_2, y_2, t_2)$ in $\R^2 \times \R$ by:
\begin{equation}
  (x_2,y_2) = (-1)^{k}(x_1, y_1)\label{eq.def_of_sim_1}
\end{equation}
and 
\begin{equation}
  t_2 = t_1 + k\label{eq.def_of_sim_2}
\end{equation}
for some $k \in \Z$.
Let $\widetilde{M} = \R^3 / \sim$ and let $q : \R^3 \to \widetilde{M}$ be its quotient map.
By \cite[Theorem 21.13]{lee2012smooth}, we see that $q : \R^3 \to \widetilde{M}$ is a smooth covering map.
In particular, the quotient map $q$ is a local diffeomorphism.
Define
\[
\widetilde{X} = q_*((V(x,y,t),1)),
\]
which is a smooth vector field on $\widetilde{M}$ (note that by the properties (1) and (2) at the beginning of Section \ref{sec.constructionofX} and the fact that $q$ is a local diffeomorphism, it is clear that $\widetilde{X}$ is well-defined).
Let $\tilde{\phi}$ be a smooth flow on $\widetilde{M}$ generated by the vector field $\widetilde{X}$.
Then, $\tilde{\phi}$ is a modified suspension flow of $\bar{f}$ in the following sense:
\begin{pro}\label{pro.tildephi_is_susp_of_barf}
  Let $i : \R^2 \to \widetilde{M}$ be $i(x,y) = q(x,y,0)$.
  Then, the following diagram commutes: 
  \begin{equation}
    \xymatrix{
    \R^2 \ar[r]^{\bar{f}} \ar[d]_{i} & \R^2 \ar[d]^{i} \\
    \widetilde{M} \ar[r]_{\tilde{\phi}(1,\cdot)} & \widetilde{M} \ar@{}[lu]|{\circlearrowright}
    }\label{eq.cumm_diag_ol_f}
  \end{equation}
  for all $(x, y) \in \R^2$.
  Equivalently, the Poincar\'e map of $\tilde{\phi}$ on the global section $q(\R^2 \times \{0\})$ is $\bar{f}$.
\end{pro}
\begin{prf}
  The vector field $(V(x,y, t), 1)$ on $\R^2 \times [0, 1]$ induces a flow such that the point $(x, y, 0)$ moves to $(f_1(x, y), 1)$ at time $1$ (see Figure \ref{figure.suspensionoff}).
  Thus, $\tilde{\phi}(1, q(x, y, 0)) = q(f(x, y), 1) = q(\bar{f}(x, y), 0)$.
  \qed
\end{prf}

Let $P : \R^3 \to \R^4$ be a smooth map defined by
\[
  P(x,y,\theta) = R_{\theta,\theta}(x,y).
\]
Then for all $(x_1,y_1,t_1), (x_2,y_2,t_2) \in \R^3$ satisfying \eqref{eq.def_of_sim_1} and \eqref{eq.def_of_sim_2},
\begin{align}
  P(x_2,y_2,t_2)
   & = R_{t_2,t_2}(x_2, y_2)                   \\
   & = R_{t_1+k,t_1+k}((-1)^k x_1, (-1)^k y_1) \\
   & = R_{t_1,t_1}(x_1, y_1)                   \\
   & = P(x_1,y_1,t_1).
\end{align}
Thus, by the universal property of the quotient maps, there is a continuous map $\widetilde{P} : \widetilde{M} \to \R^4$ such that the following diagram commutes:
\[
  \xymatrix{
  \R^3 \ar[dr]^{P} \ar[d]_{q} &  \\
  \widetilde{M}  \ar[r]_{\widetilde{P}} & \R^4
  }
\]
Since $q$ is a local diffeomorphism, $\widetilde{P}$ is differentiable.

Now, let us prove properties {\bf (P2)}, {\bf (P3)} and {\bf (P4)}.
Since $(0,0) \in \CR(\bar{f})$ is not isolated, the corresponding closed orbit $\{q(0,0,t) ; 0 \leq t \leq 1\}$ is accumulated by regular points of $\CR(\tilde{\phi})$. This and the definition of $\widetilde{P}$ imply property {\bf (P4)}.
For all $(x,y,t) \in \R^2 \times [0,1]$, we have
\begin{align}
  P_*((V(x,y,t),1))
   & = DP(x,y,t)(V(x,y,t),1)                                                  \\
   & = (V_1(x,y,t)\cos \pi t,\, V_1(x,y,t)\sin \pi t, 0,0)                      \\
   & + (0,0,V_2(x,y,t)\cos \pi t,\, V_2(x,y,t)\sin \pi t)                       \\
   & + (-\pi x\sin \pi t,\, \pi x\cos \pi t,\, -\pi y\sin \pi t,\, \pi y\cos \pi t) \\
   & = R_{t,t}(V(x,y,t)) + R_{t+1/2,t+1/2}(\pi x,\pi y)                       \\
   & = Z(R_{t,t}(x,y)) + Y(R_{t,t}(x,y))+(0,0,0,0)                                      \\
   & = X(R_{t,t}(x,y)) = X(P(x,y,t)).
\end{align}
This and the definition of $\tilde{\phi}$ imply $\widetilde{P}_* (q_*(V(x,y,t),1)) = X(\widetilde{P}(x,y,t))$, proving property {\bf (P2)}.

To prove property {\bf (P3)}, we need the following lemma:
\begin{lem}\label{lem.CR_phi_K_is_P_CR_tilde_phi}
  \[
    \CR(\phi \vert_K) = \widetilde{P}(\CR(\tilde{\phi})).
  \]
\end{lem}
\begin{prf}
  Since $\widetilde{P}$ is continuous map with $\widetilde{P}(\widetilde{M}) = K$, every pseudotrajectory of $\tilde{\phi}$ is mapped to that of $\phi \vert_K$.
  So, $\widetilde{P}(\CR(\tilde{\phi})) \subset \CR(\phi \vert_K)$.
  
  To prove the reverse inclusion, suppose, by contradiction, that there is $p \in \CR(\phi \vert_K) \setminus \widetilde{P}(\CR(\tilde{\phi}))$.
  As $\widetilde{P}$ is surjective onto $K$, there exists $q \in \widetilde{P}^{-1}(\CR(\phi \vert_K)) \setminus \CR(\tilde{\phi})$ such that $\widetilde{P}(q) = p$.
  By the choice of $\widetilde{P}$, we see that the $\tilde{\phi}$-invariant set $\widetilde{P}^{-1}(\CR(\phi \vert_K))$ is compact.
  
  From Proposition \ref{pro.tildephi_is_susp_of_barf}, it follows that $\tilde{\phi}$ is the modified suspension flow of Smale horseshoe. Therefore, every orbit of $\tilde{\phi}$ corresponds to some orbit of $\bar{f}$ via the Poincar\'e map and we have
  \[
    \CR(\tilde{\phi}) = \{ \tilde{\phi}(t, x) ; t \in [0,1], x \in i(\CR(\bar{f})) \},
  \]
  where $i$ is the map defined in Proposition \ref{pro.tildephi_is_susp_of_barf}.
  Thus, the orbit of $q$ under $\tilde{\phi}$ corresponds to a bounded orbit of $\bar{f}$ that does not intersect $CR(\bar{f})$.
  Such a bounded orbit must converge to the attracting fixed points of $\bar{f}$. Thus, the orbit of $q$ must converge to the corresponding attracting closed orbits of $\tilde{\phi}$.
  Similarly, by property {\bf (P2)}, the orbit of $p$ must converge to the attracting closed orbits of $\phi \vert_K$.
  Thus, every $d$-pseudotrajectory from $p$ with sufficiently small $d$ must converge to the attracting closed orbits, implying that $p \notin \CR(\phi \vert_K)$.
  This is a contradiction, and the proof is complete.
  \qed
\end{prf}
Now, using \cite[Theorem 1.5.36]{fisher2019hyperbolic} for $\phi$ and Proposition \ref{pro.CRphi_bdd_inK}, we have
\[
  \CR(\phi \vert_K) \subset \CR(\phi) = \CR(\phi \vert_{\CR(\phi)}) \subset \CR(\phi \vert_{K}),
\]
where the latter inclusion follows from Proposition \ref{pro.CRphi_bdd_inK}.
This and Lemma \ref{lem.CR_phi_K_is_P_CR_tilde_phi} yield property {\bf (P3)}.
\section{Proof of Theorem \ref{thm.factor}}\label{sec.prf_of_thm_factor}
Consider topological flows $(\Lambda_1, \phi_1)$ and $(\Lambda_2, \phi_2)$ on compact metric spaces $\Lambda_1$ and $\Lambda_2$.
Let $\dist_1$, $\dist_2$ be the distance functions of $\Lambda_1$ and $\Lambda_2$, respectively.
Suppose that there is a continuous surjection $P : \Lambda_1 \to \Lambda_2$ by which the following diagram commutes:
\[
  \xymatrix{
  \Lambda_1 \ar[r]^{\phi_1(t,\cdot)} \ar[d]_{P} & \Lambda_1 \ar[d]^{P} \\
  \Lambda_2 \ar[r]_{\phi_2(t,\cdot)} & \Lambda_2 \ar@{}[lu]|{\circlearrowright}
  }
\]
for all $t \in \R$.
Assume that there exists a closed orbit $\gamma \subset \Lambda_1$ of $\phi_1$ such that $P(\gamma)$ is a singularity of $\phi_2$, the restriction of $P$ to $\Lambda_1 \setminus \gamma$ is a bijection onto $\Lambda_2 \setminus P(\gamma)$, and $\phi_1$ has the standard shadowing property.
Given $\e_0 > 0$,
choose $\e_1 > 0$ so that if $\dist_1(x,y) < \e_1$ for some $x,y \in \Lambda_1$, then $\dist_2(P(x),P(y)) < \e_0/2$.
Let $T_0 > 0$ be the minimal period of $\gamma$.
Using the standard shadowing property of $\phi_1$, we may take $d_0 > 0$ so that every $d_0$-pseudotrajectory can be $\e_1$-standard shadowed.

In the following lemma, we prove that if two pseudotrajectories come close to the closed orbit $\gamma$, then we can connect these pseudotrajectories via $\gamma$.
\begin{lem}\label{lem.po_change}
  There exists $d_1 > 0$ such that if $d_1$-pseudotrajectories $\xi_1, \xi_2$ of $\phi_1$ and $p \in \gamma$ satisfy
  \[
    \dist_1(\xi_1(t_1), \gamma) < d_1, \quad
    \dist_1(\xi_2(t_2), \gamma) < d_1
  \]
  for some $t_1, t_2 \in \R$ with $t_2 - t_1 \geq T_0/d_1$, then there is a $d_0$-pseudotrajectory $\xi$ such that
  \[
    \xi(t) =
    \begin{cases}
      \xi_1(t), & t \leq t_1, \\
      \xi_2(t), & t \geq t_2
    \end{cases}
  \]
  and $\xi(t) \in \gamma$ for all $t \in (t_1, t_2)$.
\end{lem}
\begin{prf}
  There exist $p \in \gamma$ and $T \in [0, T_0)$ such that 
  \[
    \dist_1(\xi_1(t_1), p) < d_1, \quad
    \dist_1(\xi_2(t_2), \phi(T,p)) < d_1.
  \]
  There is $N > 0$ such that
  \[
    \lvert t_2 - t_1 - (NT_0 + T) \rvert \leq T_0.
  \]
  Since $t_2 - t_1 \geq T_0/d_1$, 
  \[
    \left\lvert \frac{t_2 - t_1 - (NT_0 + T)}{t_2 - t_1} \right\rvert
    \leq \frac{\lvert t_2 - t_1 - (NT_0 + T) \rvert}{T_0/d_1}
    \leq d_1.
  \]
  Thus, taking $d_1$ smaller, we may assume that the function $\xi : \R \to \Lambda_1$ defined by
  \[
    \xi(t) =
    \begin{cases}
      \xi_1(t),                                                 & t \leq t_1,       \\
      \phi_1 \left( \frac{NT_0+T}{t_2-t_1}(t - t_1), p \right), & t \in (t_1, t_2), \\
      \xi_2(t),                                                 & t \geq t_2
    \end{cases}
  \]
  is a $d_0$-pseudotrajectory.
  \qed
\end{prf}

Let
\begin{equation}
  B \subset \{ p \in \Lambda_2 ; \dist_2(p, P(\gamma)) < \e_0/2 \}\label{eq.def_of_U1}
\end{equation}
be a closed neighborhood of the singularity $P(\gamma)$ such that 
\[
  P^{-1}\left( B \right) \subset \{p \in \Lambda_1 ; \dist_1(p,\gamma) < d_1 \}
\]
(such a neighborhood does exist since the continuous image of compact set $P(\{p \in \Lambda_1 ; \dist_1(p,\gamma) \geq d_1 \})$ is closed, and disjoint from $\gamma$).
Then, choose $d_2 \in (0,d_1)$ and a neighborhood $U \subset B$ of $P(\gamma)$ such that if $d_2$-pseudotrajectory $\xi$ with $\xi(t_0) \in U$ for some $t_0 \in \R$, then
\begin{equation}
\xi(t) \in B\label{eq.def_of_U2}
\end{equation}
for all $t$ with $\lvert t - t_0 \rvert < T_0/d_1$.
Since $P$ restricted to the compact set $P^{-1}(\Lambda_2 \setminus U)$ is homeomorphism onto its image, taking $d_2$ smaller if necessary, we may assume that if a $d_2$-pseudotrajectory $\xi$ satisfies $\xi(t) \in \Lambda_2 \setminus U$ for all $t \in [t_0, t_1]$, then $P^{-1} \circ \xi$ is a $d_1$-pseudotrajectory on $[t_0, t_1]$.

Let $\xi$ be a $d_2$-pseudotrajectory of $\phi_2$.
Define
\[
  S = \{ t \in \R ; \xi(t) \in U \}
\]
and
\[
  \widetilde{S} = \{ s \in \R ; \lvert s - t \rvert < T_0/d_1 \text{ for some } t \in S \}.
\]
For $l,r \in \R \cup \{-\infty, \infty\}$, let
\[
  (l,r) = \{t \in \R ; l < t < r\}.
\]
Then $\widetilde{S}$ can be divided into a collection of open intervals $\{ (l_n, r_n) ; a_S < n < b_S, n \in \Z \}$.
Here, note that:
\begin{enumerate}
  \item The case where $a_S = -\infty$ and $b_S = \infty$ is admitted;
  \item If $a_S > -\infty$ (resp. $b_S < \infty$), then $l_{a_S+1} = -\infty$ (resp. $r_{b_S-1} = \infty$) is permissible;
  \item By the definition of $\widetilde{S}$, we have $r_n - l_n \geq T_0/d_1$ for all $n \in \Z$ with $a_S < n < b_S$.
\end{enumerate}
\begin{figure}[h]
  \centering
  \begin{minipage}[b]{1\columnwidth}
    \centering
    \includegraphics[width=12cm]{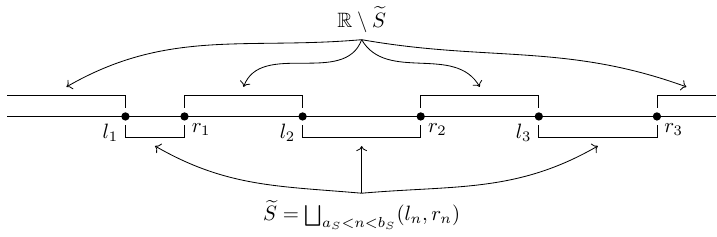}
    \subcaption{Bounded $\widetilde{S}$.}\label{figure.T1}
    \vfill
    \centering
    \includegraphics[width=12cm]{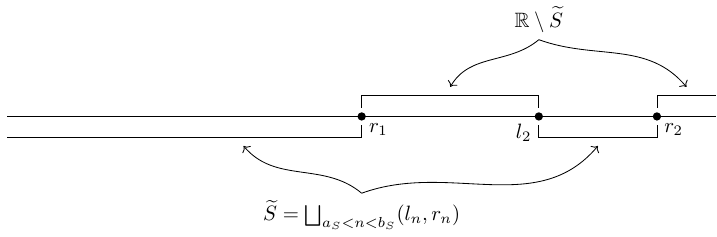}
    \subcaption{Unbounded $\widetilde{S}$.}\label{figure.T2}
  \end{minipage}
\end{figure}
Let us construct a $d_0$-pseudotrajectory $\tilde{\xi} : \R \to \Lambda_1$ of $\phi_1$.
Define $\xi' : \R \setminus \widetilde{S} \to \Lambda_1$ by
\[
  \xi'(t) = P^{-1} \circ \xi(t).
\]
Then $\xi'$ is thought to be a union of $d_1$-pseudotrajectories of closed intervals.
Let $n \in \Z$ with $a_S < n < b_S$.
By the choice of $U$ and $\widetilde{S}$, we have $\xi(l_n), \xi(r_n) \in B$.
Thus, there exists $p_n, q_n \in \gamma$ such that $\dist_1(P^{-1} \circ \xi(l_n), p_n), \dist_1(P^{-1} \circ \xi(r_n), q_n) < d_1$.
Applying Lemma \ref{lem.po_change} to $\xi'$ for each $n$, we obtain that there is a $d_0$-pseudotrajectory $\tilde{\xi}$ such that
\[
  \tilde{\xi}(t) \in \gamma, \quad t \in (l_n, r_n) \text{ for some $a_S < n < b_S$}
\]
and
\[
  \tilde{\xi}(t) = P^{-1} \circ \xi(t)
\]
for all $t \in \R \setminus \widetilde{S}$.
Then by the choice of $d_0$, there exist $x \in \Lambda_1$ and $h \in \Rep(\e_1)$ with $h(0) = 0$ such that
\begin{equation}
  \dist_1(\tilde{\xi}(t), \phi_1(h(t), x)) < \e_1, \quad t \in \R.\label{eq.shadowing_for_tilde_xi}
\end{equation}
In order to prove Theorem \ref{thm.factor}, it is enough to prove
\begin{equation}
  \dist_2(\xi(t), \phi_2(h(t), P(x))) < \e_0, \quad t \in \R.\label{eq.shadowing_for_xi}
\end{equation}

By \eqref{eq.def_of_U2}, we have $\xi(t) \in B$ for all $t \in [l_n, r_n]$.
This and \eqref{eq.def_of_U1} imply
\[
  \dist_2(\xi(t), P(\gamma)) < \frac{\e_0}{2}, \quad t \in [l_n, r_n].
\]
It follows from \eqref{eq.shadowing_for_tilde_xi} and the choice of $\e_1$ that
\[
  \dist_2(P(\gamma), \phi_2(h(t), P(x)))
  = \dist_2(P(\tilde{\xi}(t)), P(\phi_1(h(t), x))) < \frac{\e_0}{2}, \quad t \in [l_n, r_n].
\]
Combining these two inequalities above, we obtain \eqref{eq.shadowing_for_xi} for all $t \in [l_n, r_n]$ with $a_S < n < b_S$.

It remains to prove that \eqref{eq.shadowing_for_xi} holds for all $t \in \R \setminus \widetilde{S}$.
By \eqref{eq.shadowing_for_tilde_xi} and the choice of $\e_1$,
\[
  \dist_2(\xi(t), \phi_2(h(t), P(x)))
  = \dist_2(P(\tilde{\xi}(t)), P(\phi_1(h(t), x))) < \frac{\e_0}{2}
\]
for all $t \in \R \setminus \widetilde{S}$, finishing the proof of Theorem \ref{thm.factor}.

\section*{acknowledgements}
The author is grateful to my advisor S. Hayashi for his constructive suggestions and continuous support.
This work was supported by JSPS KAKENHI Grant Number JP23KJ0657.

\bibliography{math}

\providecommand{\bysame}{\leavevmode\hbox to3em{\hrulefill}\thinspace}
\providecommand{\MR}{\relax\ifhmode\unskip\space\fi MR }
\providecommand{\MRhref}[2]{%
  \href{http://www.ams.org/mathscinet-getitem?mr=#1}{#2}
}
\providecommand{\href}[2]{#2}
\begin{thebibliography}{10}

\bibitem{ArbietoAttachedSing}
Alexander Arbieto, Andr{\'e}s M.~L{\'o}pez, Elias Rego, and Yeison S{\'a}nchez,
  \emph{On the shadowableness of flows with hyperbolic singularities},
  Mathematische Annalen \textbf{390} (2024), no.~1, 417--437.

\bibitem{fisher2019hyperbolic}
Todd Fisher and Boris Hasselblatt, \emph{Hyperbolic flows}, 2019.

\bibitem{KomuroLorenz}
Motomasa Komuro, \emph{Lorenz attractors do not have the pseudo-orbit tracing
  property}, Journal of the Mathematical Society of Japan \textbf{37} (1985),
  no.~3, 489--514.

\bibitem{lee2012smooth}
John~M. Lee, \emph{Introduction to smooth manifolds}, Springer, 2012.

\bibitem{Mura}
Sogo Murakami, \emph{Oriented and standard shadowing properties for topological
  flows}, Tokyo Journal of Mathematics (2023).

\bibitem{P.S.T.}
Sergei~Yu. Pilyugin and Sergey Tikhomirov, \emph{Lipschitz shadowing implies
  structural stability}, Nonlinearity \textbf{23} (2010), no.~10, 2509.

\bibitem{TomasSuspensionShadowing}
Romeo~F. Thomas, \emph{Stability properties of one-parameter flows},
  Proceedings of the London Mathematical Society \textbf{s3-45} (1982), no.~3,
  479--505.

\bibitem{WANG2014210}
Qiqi Wang, Rui Hu, and Patrick Blonigan, \emph{Least squares shadowing
  sensitivity analysis of chaotic limit cycle oscillations}, Journal of
  Computational Physics \textbf{267} (2014), 210--224.

\bibitem{Wu:2018wo}
Xinxing Wu, Lidong Wang, and Jianhua Liang, \emph{The chain properties and
  average shadowing property of iterated function systems}, Qualitative Theory
  of Dynamical Systems \textbf{17} (2018), no.~1, 219--227.

\bibitem{WenWenNoShadowing}
Lan~Wen Xiao~Wen, \emph{No-shadowing for singular hyperbolic sets with a
  singularity.}, Discrete and Continuous Dynamical Systems \textbf{40(10):}
  (2020), 6043--6059.

\end{thebibliography}
\bibliographystyle{amsplain}
\end{document}